\documentclass[12pt]{article}
\usepackage{amssymb}
\usepackage{mathrsfs}
\usepackage{pst-poly}     
\usepackage{cite}

\newtheorem{thm}{Theorem}[section]
\newtheorem{lem}[thm]{Lemma}
\newtheorem{Def}[thm]{Definition}

\newtheorem{prop}[thm]{Proposition}

\newenvironment{pf}[1][Proof]{\noindent\textbf{#1.} }{\hfill\rule{1mm}{2mm}}

\makeatletter \@addtoreset{equation}{section} \makeatother

\begin{document}
\title{\bf Fault-tolerant analysis of augmented cubes\thanks {The work was supported by NNSF of China (No.11071233 and No.11071223) and ZSDZZZZXK08.}}
\author
{Meijie Ma\\
{\small Department of Mathematics}\\
{\small  Zhejiang Normal University}\\
{\small Jinhua, 321004, China} \\ \\
Yaxing Song  \quad Jun-Ming Xu\footnote{Corresponding author:
xujm@ustc.edu.cn}\\
{\small School of Mathematical Sciences}\\
{\small University of Science and Technology of China}\\
{\small Hefei, 230026, China}  \\
 }
\date{}
 \maketitle

\begin{abstract}
The augmented cube $AQ_n$, proposed by Choudum and Sunitha [S. A.
Choudum, V. Sunitha, Augmented cubes, Networks 40 (2) (2002) 71-84],
is a $(2n-1)$-regular $(2n-1)$-connected graph $(n\ge 4)$. This
paper determines that the 2-extra connectivity of $AQ_n$ is $6n-17$
for $n\geq 9$ and the 2-extra edge-connectivity is $6n-9$ for $n\geq
4$. That is, for $n\geq 9$ (respectively, $n\geq 4$), at least
$6n-17$ vertices (respectively, $6n-9$ edges) of $AQ_n$ have to be
removed to get a disconnected graph that contains no isolated
vertices and isolated edges. When the augmented cube is used to
model the topological structure of a large-scale parallel processing
system, these results can provide more accurate measurements for
reliability and fault tolerance of the system.

\end{abstract}

{\noindent{\bf Keywords:} Combinatorics, Fault-tolerant analysis,
Augmented cube, Extra connectivity, Extra edge-connectivity

\section{Introduction}
It is well known that the underlying topology of an interconnection
network can be modeled by a graph $G=(V, E)$, where $V$ is the set
of processors and $E$ is the set of communication links in the
network. For all the graph terminologies and notations not defined
here, we follow\cite{q2}. Then we use graphs and networks
interchangeably in this paper.

A set of vertices (respectively, edges) $S$ of $G$ is called a
vertex-cut (respectively, an edge-cut) if $G-S$ is disconnected. The
connectivity $\kappa(G)$ (respectively, the edge-connectivity
$\lambda(G)$) of $G$ is defined as the minimum cardinality of a
vertex-cut (respectively, an edge-cut) $S$. And it is known to all
that the connectivity $\kappa(G)$ and the edge-connectivity
$\lambda(G)$ are two important parameters to measure reliability and
fault tolerance of the network. These parameters, however, have an
obvious deficiency, that is, they tacitly assume that all elements
in any subset of $G$ can potentially fail at the same time, which
happens almost impossible in practiced applications of networks. In
other words, in the definition of $\kappa(G)$ and $\lambda(G)$,
absolutely no conditions or restrictions are imposed either on the
set $S$ or on the components of $G-S$. Consequently, to compensate
for these shortcomings, it seems natural to generalize the classical
connectivity by adding some conditions or restrictions on the set
$S$ and the components of $G-S$.

In \cite{es89,q5}, Esfahanian and Hakimi generalized the notion of
connectivity by suggesting the concept of restricted connectivity in
point of view of network applications. A set $S\subset V(G)$
(respectively, $F\subset E(G)$) is called a restricted vertex-set
(respectively, edge-set) if it does not contain the neighbor-set of
any vertex in $G$ as its subset. A restricted vertex-set $S$
(respectively, edge-set $F$) is called a restricted vertex-cut
(respectively, edge-cut) if $G-S$ is disconnected. The restricted
connectivity $\kappa'(G)$ (respectively, edge-connectivity
$\lambda'(G)$) is the minimum cardinality of a restricted vertex-cut
(respectively, edge-cut) in $G$, if any, and does not exist
otherwise, denoted by $+\infty$.

However, the maximum difficult for computing the restricted
connectivity of a graph $G$ is to check that a vertex-cut does not
contain the neighbor-set of any vertex in $G$ as its subset. Thus,
only a little knowledge of results has been known on $\kappa'(G)$ or
$\lambda'(G)$ even for particular classes of graphs. For example, Xu
and Xu~\cite{q13} studied the $\lambda'(G)$ for a vertex-transitive
graph $G$, Esfahanian~\cite{es89} determined
$\kappa'(Q_n)=\lambda'(Q_n)=2n-2$ for the hypercube $Q_n$ and $n\ge
3$.

To avoid this difficult, one slightly modified the concept of a
restricted vertex-set $S$ by replacing the term ``any vertex in $G$"
in the condition by  the term ``any vertex in $G-S$". We call in
this sense the connectivity as the super connectivity, denoted by
$\kappa_1(G)$ for the super connectivity and $\lambda_1(G)$ for the
super edge-connectivity (see, for example, \cite{q3, bbst81, q11,
q12, xlmh05}). Clearly, $\kappa_1(G)\leq \kappa'(G)$ and
$\lambda_1(G)=\lambda'(G)$ if they exist. The super connectivity of
some graphs determined by several authors. For example, for the
hypercube $Q_n$, the twisted cube $TQ_n$, the cross cube $CQ_n$, the
M\"obius cube $MQ_n$ and the locally twisted cube $LTQ_n$, Xu et
al~\cite{xww10} showed that their super connectivity and the super
edge-connectivity are all $2n-2$; for the star graph $S_n$, Hu and
Yang~\cite{hy97} proved that $\kappa_1(S_n)=2n-4$ for $n\geq 3$; for
the augmented cube $AQ_n$, Ma, Liu and Xu~\cite{mlx08, mtxl09}
determined $\kappa_1(AQ_n)=4n-8$ for $n\geq6$ and
$\lambda_1(AQ_n)=4n-4$ for $n\geq 2$; for the $(n,k)$-star graphs
$S_{n,k}$, Yang et al.~\cite{ylg10} proved that
$\kappa_1(S_{n,k})=n+k-3$; for the $n$-dimensional alternating group
graph $AG_n$, Cheng et al.~\cite{cls10} determined
$\kappa_1(AG_n)=4n-11$ for $n\geq 5$.

Observing that every component of $G-S$ contains at least two
vertices if $S$ is a restricted vertex-set of $G$, F\`{a}brega and
Fiol \cite{q7} introduced the $h$-extra connectivity of $G$. A
vertex-cut (respectively, an edge-cut) $S$ of G is called an
$h$-vertex-cut (respectively, an $h$-edge-cut) if every component of
$G-S$ has more than $h$ vertices. The $h$-extra connectivity
$\kappa_h(G)$ (respectively, $h$-extra edge-connectivity
$\lambda_h(G)$) defined as $\min\{|S|$: S is an $h$-vertex-cut
(respectively, $h$-edge-cut) of $G\}$. Clearly,
$\kappa_0(G)=\kappa(G)$ and $\lambda_0(G)=\lambda(G)$ for any graph
$G$ if $G$ is not a complete graph. Thus, as a new measurement for
reliability and fault tolerance of the large-scale parallel
processing system, the $h$-extra connectivity is more accurate than
the classical connectivity and has recently received much attention.
For example, for the hypercube $Q_n$, Xu et al~\cite{xzhz05, zx06}
determined $\kappa_2(Q_n)=3n-5$ and $\lambda_2(Q_n)=3n-4$ for $n\geq
4$; for the folded hypercube $FQ_n$, Zhu et al~\cite{zxhx07}
determined $\kappa_2(FQ_n)=3n-2$ for $n\geq 8$ and
$\lambda_2(FQ_n)=3n-1$ for $n\geq 5$; for the star graph $S_n$, Wan,
Zhang~\cite{wz09} determined $\kappa_2(S_n)=6(n-3)$ for $n\geq 4$;
for the $(n,k)$-star graphs $S_{n,k}$, Yang et al.~\cite{ylg10}
proved that $\kappa_2(S_{n,k})=n+2k-5$ for $2\leq k\leq n-2$; Zhang
et al.~\cite{zxy10} proved $\kappa_2(AG_n)=6n-18$ for $n\geq 5$.


In this paper, we study the augmented cube $AQ_n$ and determine
$\kappa_2(AQ_n)=6n-17$ for $n\geq 9$ and $\lambda_2(AQ_n)=6n-9$ for
$n\geq 4$.

The rest of this paper is organized as follows. In Section 2,
we recall the structure of $AQ_n$, and some definitions and
lemmas. The main results are given in Section 3.
Finally, we conclude our paper in Section 4.


\section{Definitions and lemmas}

Let $n$ be a positive integer. The {\it $n$-dimensional augmented
cube}, denoted by  $AQ_n$, proposed by Choudum and
Sunitha~\cite{cs00, cs01, cs02}, having $2^{n}$ vertices, each
labeled by an $n$-bit binary string, that is,
$V(AQ_n)=\{x_{n}x_{n-1}\cdots x_{1}:\ x_{i} \in\{0,1\},\ 1\leq i\leq
n\}$, can be defined recursively as follows.

\begin{Def}\label{Def2.1}\
$AQ_1$ is a complete graph $K_2$ with the vertex set $\{0, 1\}$. For
$n\geq 2$, $AQ_n$ is obtained by taking two copies of the augmented
cube $AQ_{n-1}$, denoted by $AQ_{n-1}^0$ and $AQ_{n-1}^1$, where
$V(AQ_{n-1}^0)= \{0x_{n-1}\ldots x_2x_1: x_i\in\{0,1\},\ 1\leq i\leq
n-1\}$ and $V(AQ_{n-1}^1)= \{1x_{n-1}\ldots x_2x_1: x_i\in\{0,1\},\
1\leq i\leq n-1\}$, and a vertex $X=0x_{n-1}\ldots x_2x_1$ of
$AQ_{n-1}^0$ being joined to a vertex $Y=1y_{n-1}\ldots y_2y_1$ of
$AQ_{n-1}^1$ if and only if either

\begin{enumerate}
\item[(i)] $x_i=y_i$ for $1\leq i\leq n-1$,
or

\item[(ii)] $x_i=\bar y_i$ for $1\leq i\leq n-1$.
\end{enumerate}
\end{Def}

The graphs shown in Figure \ref{f1} are the augmented cubes $AQ_1$,
$AQ_2$ and $AQ_3$, respectively.

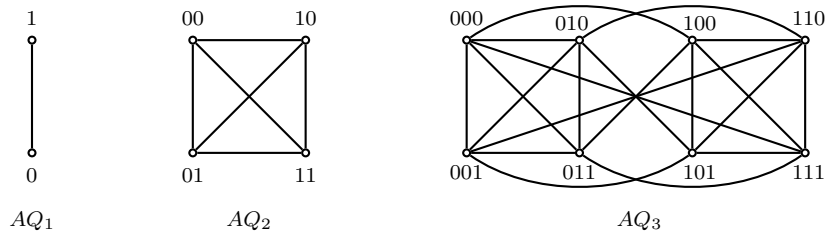
\begin{figure}[h]

\begin{pspicture}(-1.5,0)(.5,3.5)
\psset{radius=.06}

\Cnode(1,1){0}\rput(1,.7){\scriptsize0}
\Cnode(1,2.5){1}\rput(1,2.8){\scriptsize1} \ncline{1}{0}
\rput(1,0.1){\scriptsize $AQ_1$}
\end{pspicture}
\begin{pspicture}(-1.5,0)(3,3.5)
\psset{radius=.06}

\Cnode(1,1){01}\rput(1,.7){\scriptsize01}
\Cnode(2.5,1){11}\rput(2.5,.7){\scriptsize11}
\Cnode(2.5,2.5){10}\rput(2.5,2.8){\scriptsize10}
\Cnode(1,2.5){00}\rput(1,2.8){\scriptsize00}

\ncline{00}{01}\ncline{00}{10}\ncline{10}{11}\ncline{11}{00}
\ncline{10}{01}\ncline{11}{01} \rput(1.75,0.1){\scriptsize $AQ_2$}
\end{pspicture}
\begin{pspicture}(-.5,0)(5.5,3.5)
\psset{radius=.06}

\Cnode(1,1){001}\rput(1,.75){\scriptsize001}
\Cnode(2.5,1){011}\rput(2.5,.75){\scriptsize011}
\Cnode(2.5,2.5){010}\rput(2.4,2.73){\scriptsize010}
\Cnode(1,2.5){000}\rput(1,2.8){\scriptsize000}

\ncline{000}{001}\ncline{000}{010}\ncline{010}{011}\ncline{011}{000}
\ncline{010}{001}\ncline{011}{001}

\Cnode(4,1){101}\rput(4.1,.75){\scriptsize101}
\Cnode(5.5,1){111}\rput(5.55,.75){\scriptsize111}
\Cnode(5.5,2.5){110}\rput(5.55,2.8){\scriptsize110}
\Cnode(4,2.5){100}\rput(4.1,2.73){\scriptsize100}

\ncline{100}{101}\ncline{100}{110}\ncline{110}{111}\ncline{111}{100}
\ncline{110}{101}\ncline{111}{101}
\ncline{010}{101}\ncline{100}{011}\ncline{000}{111}\ncline{001}{110}
\nccurve[angleA=35,angleB=145]{000}{100}
\nccurve[angleA=35,angleB=145]{010}{110}
\nccurve[angleA=-35,angleB=-145]{001}{101}
\nccurve[angleA=-35,angleB=-145]{011}{111}

\rput(3.3,0.1){\scriptsize $AQ_3$}

\end{pspicture}

\caption{\label{f1}\footnotesize { Three augmented cubes $AQ_1$,
$AQ_2$ and $AQ_3$}}
\end{figure}

For convenience, we can express the recursive structure of $AQ_n$ as
$AQ_n=L\odot R$, where $L=AQ^{0}_{n-1}$ and $R=AQ^{1}_{n-1}$. Then
we call the edges between $L$ and $R$ {\it crossed edges}. Obviously
every vertex is incident to exactly two
crossed edges. Let $X=x_{n}x_{n-1}\cdots x_{1}$ be an $n$-bit binary
string. And for $1\leq i\leq n$, let
 $$
 \begin{array}{rl}
 & X_i=x_{n}x_{n-1}\cdots x_{i+1}\bar{x}_{i}x_{i-1}\cdots x_{1},\\
 &\overline{X}_i=x_{n}x_{n-1}\cdots x_{i+1}\bar{x}_{i}\bar{x}_{i-1}\cdots \bar{x}_{1}.
 \end{array}
 $$
Obviously, $X_1=\overline{X}_1$,
$(X_i)_i=X=\overline{(\overline{X}_i)}_{i}$. According to
Definition~\ref{Def2.1}, we can directly obtain a useful
characterization of adjacency.

\begin{prop}\label{prop2.2}\ Assume that $X=x_{n}x_{n-1}\cdots
x_{1}$ and $Y=y_{n}y_{n-1}\cdots y_{1}$ are two distinct vertices in
$AQ_n$. Then $X$ and $Y$ are adjacent if and only if either

i)\ there exists an integer $i\ (1\leq i\leq n)$ such that $Y=X_i$,
or

ii)\ there exists an integer $i\ (2\leq i\leq n)$ such that
$Y=\overline{X}_i$.
\end{prop}

By Proposition~\ref{prop2.2}, an alternative definition of $AQ_n$
can be stated as follows.

\vskip6pt\begin{Def}\label{Def2.2}\ The augmented cube $AQ_n$ of
dimension $n$ has $2^n$ vertices. Each vertex is labeled by a unique
$n$-bit binary string as its address. Two vertices $X$ and $Y$ are
joined if and only if either

\begin{enumerate}
\item[(i)] there exists an integer
$i$ with $1\leq i\leq n$ such that $Y=X_i$; in this case, the edge
is called a hypercube edge of dimension $i$, denoted by $XX_i$, or

\item[(ii)] there exists an integer
$i$ with $2\leq i\leq n$ such that $Y=\overline X_i$; in this case,
the edge is called a complement edge of dimension $i$, denoted by
$X\overline X_i$.

\end{enumerate}
\end{Def}

\begin{lem}\label{lem2.3} \textnormal{(Choudum and
Sunitha~\cite{cs02})}\ $AQ_n$ is $(2n-1)$-regular
$(2n-1)$-connected for $n\geq4$, however, $\kappa(AQ_3)=4$ for $n=3$.
\end{lem}

\begin{lem}\label{lem2.5a} \textnormal{(Ma, Liu and Xu~\cite{mlx08, mtxl09})}\
$\kappa_1(AQ_n)=4n-8$ for $n\geq6$ and $\lambda_1(AQ_n)=4n-4$ for
$n\geq2$.
\end{lem}

\begin{lem}\label{lem2.4}\
Any two adjacent vertices in $AQ_n$ have either two or four common
neighbors for $n\geq3$.
\end{lem}

\begin{pf}\ Let $X$ and $Y$
be two adjacent vertices in $AQ_n$. Then $Y$ is either $X_i$ or
$\overline{X}_i$ by Proposition~\ref{prop2.2}. If $Y=X_i$ for some
$i$ with $1\leq i\leq n$, that is, $XX_i$ is a hypercube edge of
dimension $i$, then we have
 \begin{equation}\label{e2.1}
  N_{AQ_n}(X)\cap N_{AQ_n}(X_i)=\left\{
  \begin{array}{lll}
     \{X_2,\overline{X}_2\} &\ {\rm if }\  i=1\\
      \{\overline{X}_i,\overline{X}_{i-1}\} &\ {\rm if }\ i>1,\
    \end{array} \right.
 \end{equation}
that is, $X$ and $X_i$ have exactly two common neighbors in $AQ_n$.
If $Y=\overline{X}_i$ for some $i$ with $2\leq i\leq n$, that is,
$X\overline{X}_i$ is a complement edge of dimension $i$, then we
have
 \begin{equation}\label{e2.2}
  N_{AQ_n}(X)\cap N_{AQ_n}(\overline{X}_i)=\left\{ \begin{array}{lll}
     \{X_i,X_{i+1},\overline{X}_{i-1},\overline{X}_{i+1}\} &\ {\rm if }\ 2\leq i \leq{n-1}\\
      \{\overline{X}_{n-1},X_n\}&\ {\rm if }\ i=n.\
    \end{array} \right.
 \end{equation}
In this case, $X$ and $\overline{X}_i$ have four common neighbors
for $2\leq i\leq n-1$ while have two common neighbors for $i=n$.
\end{pf}

\begin{lem}\label{lem2.5} \textnormal{(Ma, Liu and Xu~\cite{mlx08})}\
Any two vertices in $AQ_n$ have at most four common neighbors for
$n\geq3$.
\end{lem}

Let $N_{AQ_n}(T)=\cup_{U\in V(T)}N_{AQ_n}(U)\setminus V(T)$ and $E_{AQ_n}(T)=\{XY\mid XY\in E(AQ_n)\,\, \mbox {and}\,\, X\in
V(T), Y\in V(AQ_n)\setminus V(T)\}$ for any subgraph $T$ of $AQ_n$, we have the following consequences.

\begin{lem}\label{lem2.6} \
Let $P=(Y, X, Z)$ be a path of length two in $AQ_n$ between $Y$ and $Z$ for $n\geq5$.
Then $|N_{AQ_n}(P)|\geq{6n-17}$ and $|N_{AQ_n} (X) \cap N_{AQ_n} (Y)
\cap N_{AQ_n}(Z)|\leq 1$. Furthermore, if $Z=\overline{X}_n$, we
have $|N_{AQ_n}(P)|\geq{6n-15}$.
\end{lem}

\begin{pf}
According to
Proposition~\ref{prop2.2}, $Y$ and $Z$ are in $N_{AQ_n}(X)=A\cup B$,
where $A=\{X_i \mid 1\leq i\leq n\}$ and $B=\{\overline{X}_j \mid
2\leq j\leq n\}$. It is clear that $A\cap B=\emptyset$. Consider the
following three cases.

\par
\textbf{Case 1.} $\{Y, Z\}\subset A$.

In this case, $XY$ and $XZ$ are all hypercube edges of some
dimensions $i$ and $j$, respectively. Without loss of generality, we
may assume that $Y=X_i$, $Z=X_j$ for $1\leq i<j\leq n$. By
Definition~\ref{Def2.2} and Lemma~\ref{lem2.4}, we have
 $$
 \begin{array}{rl}
  & N_{AQ_n}(X)\cap N_{AQ_n}(X_i)=\left\{ \begin{array}{lll}
     \{X_2,\overline{X}_2\} &\ {\rm if }\  i=1\\
      \{\overline{X}_i,\overline{X}_{i-1}\} &\ {\rm if }\  i>1,\
    \end{array} \right.\\
 & N_{AQ_n}(X)\cap N_{AQ_n}(X_j)=\{\overline{X}_j,\overline{X}_{j-1}\}
 \end{array}
 $$
and
 $$
  N_{AQ_n}(X_i)\cap N_{AQ_n}(X_j)=\left\{ \begin{array}{lll}
     \{X,\overline{X}_2,(X_1)_{3},\overline{(X_1)}_{3}\} &\ {\rm if }\   i=1,j=2,3\\
      \{\overline{X}_i,X,(X_i)_{i+1},\overline{(X_i)}_{i+1}\} &\ {\rm if }\  i>1,j=i+1\\
      \{X,(X_i)_j\} &\ {\rm otherwise.}
    \end{array} \right.
 $$
Since $AQ_n$ is $(2n-1)$-regular and $(2n-1)$-connected for
$n\geq5$, it is not difficult to check that
\begin{equation}\label{e2.3}
|N_{AQ_n}(P)|=\left\{ \begin{array}{lll}
     6n-13 &\ {\rm if }\    i=1,j=2,3\\
      6n-13 &\ {\rm if }\   i>1,j=i+1\\
      6n-12 &\ {\rm otherwise.}
    \end{array} \right.
\end{equation}

\textbf{Case 2.} $\{Y, Z\}\subset B$.

In this case, $XY$ and $XZ$ are all complement edges of some
dimensions $i$ and $j$, respectively. Without loss of generality,
assume that $Y=\overline{X}_i$ and $Z=\overline{X}_j$ for $2\leq
i<j\leq n$, then we have
 $$
 \begin{array}{rl}
  &N_{AQ_n}(X)\cap
 N_{AQ_n}(\overline{X}_i)=\{X_i,X_{i+1},\overline{X}_{i-1},\overline{X}_{i+1}\},\\
 & N_{AQ_n}(X)\cap N_{AQ_n}(\overline{X}_j)=\left\{ \begin{array}{lll}
     \{X_j,X_{j+1},\overline{X}_{j-1},\overline{X}_{j+1}\} &\ {\rm if }\     j<n\\
      \{\overline{X}_{n-1},X_n\} &\ {\rm if }\    j=n\
    \end{array} \right.
    \end{array}
 $$
and
 $$
  N_{AQ_n}(\overline{X}_i)\cap N_{AQ_n}(\overline{X}_j)=\left\{ \begin{array}{lll}
     \{X,X_{i+1}\} &\ {\rm if }\     j=i+1,j<n\\
      \{X,X_n\} &\ {\rm if }\    j=i+1,j=n\\
      \{X,\overline{X}_{i+1},(\overline{X}_i)_{i+2},\overline{(\overline{X}_i)}_{i+2}\} &\ {\rm if }\    j=i+2,j<n\\
      \{X,\overline{X}_{n-1},(\overline{X}_{n-2})_{n},\overline{(\overline{X}_{n-2})}_{n}\} &\ {\rm if }\   j=i+2,j=n\\
      \{X,\overline{(\overline{X}_{i})}_{j}\} &\ {\rm if }\   j\geq i+3,j<n\\
      \{X,\overline{(\overline{X}_{i})}_{n}\} &\ {\rm if }\   j\geq i+3,j=n\
\end{array} \right.
$$
It is easy to compute that
\begin{equation}\label{e2.4}
|N_{AQ_n}(P)|=\left\{ \begin{array}{lll}
6n-15 &{\rm  if }j=i+1,j<n\\
6n-13 &{\rm  if }j=i+1,j=n\\
6n-17 &{\rm  if }j=i+2,j<n\\
6n-15 &{\rm  if }j=i+2,j=n\\
6n-16 &{\rm  if }j\geq i+3,j<n\\
6n-14 &{\rm  if }j\geq i+3,j=n.
\end{array} \right.
\end{equation}

\textbf{Case 3.} $Y\in A$ and $Z\in B$.

In this case, $XY$ is a hypercube edge of some dimension $i$ and
$XZ$ is a complement hypercube edge of some dimension $j$. Assume
that $Y=X_i$, $Z=\overline{X}_j$, for $1\leq i\leq n, 2\leq j\leq
n,$ we have
 $$
 \begin{array}{rl}
  & N_{AQ_n}(X)\cap N_{AQ_n}(X_i)=\left\{ \begin{array}{lll}
     \{X_2,\overline{X}_2\} &\ {\rm if }\     i=1\\
      \{\overline{X}_i,\overline{X}_{i-1}\} &\ {\rm if }\   i>1,\
    \end{array} \right.\\
  &N_{AQ_n}(X)\cap N_{AQ_n}(\overline{X}_j)=\left\{ \begin{array}{lll}
     \{X_j,X_{j+1},\overline{X}_{j-1},\overline{X}_{j+1}\} &\ {\rm if }\     j<n\\
      \{\overline{X}_{n-1},X_n\} &\ {\rm if }\    j=n\
    \end{array} \right.
    \end{array}
 $$
and
 $$
 N_{AQ_n}(X_i)\cap N_{AQ_n}(\overline{X}_j)= \\ \left\{ \begin{array}{lll}
   \{X,X_2\} &\ {\rm if }\    i=1,j=2\\
   \{X,\overline{X}_2,(X_1)_3,\overline{(X_1)}_{3}\} &\ {\rm if }\    i=1,j=3 \\
   \{X,\overline{(X_1)}_{j}\} &\ {\rm if }\    i=1,4\leq j<n\\
   \{X,\overline{(X_1)}_{n}\} &\ {\rm if }\    i=1,j=n\\
   \{X,X_1\} &\ {\rm if }\    i=j=2\\
   \{X,\overline{X}_{i-1},(X_i)_{i-1},(\overline{X}_{i})_{i-1}\} &\ {\rm if }\    3\leq i=j\leq{n-1}\\
   \{X,\overline{X}_{n-1},(X_n)_{n-1},(\overline{X}_{n})_{n-1}\} &\ {\rm if }\    i=j=n\\
   \{X,\overline{X}_{i},(X_i)_{i+1},\overline{(X_i)}_{i+1}\} &\ {\rm if }\     j=i-1,3\leq i\leq{n-1}\\
   & \quad {\rm or} \ j=i+1,2\leq i\leq{n-2}\\
   \{X,\overline{X}_n\} &\ {\rm if }\    j=i-1,i=n\\
   \{X,\overline{X}_{n-1},(X_{n-1})_{n},\overline{(X_{n-1})}_{n}\} &\ {\rm if }\     j=n,i=n-1\\
   \{X,\overline{X}_{i-1},(X_{i})_{i-1},\overline{(X_i)}_{i-2}\} &\ {\rm if }\     j=i-2,i\geq4\\
   \{X,\overline{(X_i)}_{j}\} &\ {\rm if }\   j\leq{i-3},i\geq5 \quad \\
   &\quad {\rm or} \ j\geq{i+2},j<n\\
   \{X,\overline{(X_i)}_{n}\} &\ {\rm if }\   j\geq{i+2},j=n.
\end{array} \right.
 $$
We can compute that
 \begin{equation}\label{e2.5}
|N_{AQ_n}(P)|=\left\{ \begin{array}{lll}
    6n-13  &\ {\rm if }\    i=1,j=2\\
    6n-15  &\ {\rm if }\    i=1,j=3 \\
    6n-14  &\ {\rm if }\    i=1,4 \leq j < n\\
    6n-12  &\ {\rm if }\    i=1,j=n\\
    6n-13  &\ {\rm if }\    i=j=2\\
    6n-15  &\ {\rm if }\    3 \leq i=j \leq {n-1}\\
    6n-13  &\ {\rm if }\    i=j=n\\
    6n-15  &\ {\rm if }\     j=i-1,3 \leq i \leq {n-1} \quad \\
    &\quad {\rm or} \ j=i+1,2 \leq i \leq {n-2}\\
    6n-13  &\ {\rm if }\    j=i-1,i=n\\
    6n-13  &\ {\rm if }\     j=n,i=n-1\\
    6n-15  &\ {\rm if }\     j=i-2,i \geq 4\\
    6n-14  &\ {\rm if }\    j \leq {i-3},i \geq 5 \quad \\
    &\quad  {\rm or} \ j \geq {i+2},j<n\\
    6n-12  &\ {\rm if }\    j \geq {i+2},j=n.
\end{array} \right.
\end{equation}\
in view of (\ref{e2.3}),(\ref{e2.4}) and (\ref{e2.5}), we derive
that $|N_{AQ_n}(P)|\geq{6n-17}$.

From the above, we can easily check that
$|N_{AQ_n} (X) \cap N_{AQ_n} (Y) \cap
N_{AQ_n}(Z)|\leq 1$ in all these cases.
By (\ref{e2.4}) and (\ref{e2.5}), we have $|N_{AQ_n}(P)|\geq{6n-15}$
if $Z=\overline{X}_n$.

The lemma follows.
\end{pf}

\vskip6pt\begin{lem}\label{lem2.7} \ 
Let $P=(Y, X, Z)$ be a path of length two in $AQ_n$ for $n\geq5$.
Assume that $U\in N_{AQ_n}(P)$, then
$|N_{AQ_n}(U,X,Y,Z)| \geq 8n-31$. If $Z=\overline{X}_n$, we have
$|N_{AQ_n}(U,X,Y,Z)|\geq 8n-29$.
\end{lem}

\begin{pf}\ 
By Lemma~\ref{lem2.6}, we have $|N_{AQ_n}(P)|\geq {6n-17}$, and
$|N_{AQ_n}(P)|\geq{6n-15}$ if $Z=\overline{X}_n$. Since $U$ is a
vertex in $N_{AQ_n}(P)$, $U$ has at most three neighbors in $P$

If $U$ has exactly one neighbor in $P$, by Lemma~\ref{lem2.5} we
have $|N_{AQ_n} (U)\cap N_{AQ_n}(X,Y,Z)| \leq 11$. In this case we
can easily compute that $|N_{AQ_n}(U,X,Y,Z)| \geq
6n-17-1+2n-1-1-11=8n-31$, and $|N_{AQ_n}(U,X,Y,Z)|\geq
6n-15-1+2n-1-1-11=8n-29$ if $Z=\overline{X}_n$.

If $U$ has exactly two neighbors in $P$, by Lemma~\ref{lem2.5}, we
have $|N_{AQ_n} (U)\cap N_{AQ_n}(X,Y,Z)| \leq 10$. In this case we
arrive at $|N_{AQ_n}(U,X,Y,Z)| \geq 6n-17-1+2n-1-2-10=8n-31$, and
$|N_{AQ_n}(U,X,Y,Z)|\geq 6n-15-1+2n-1-2-10=8n-29$ if
$Z=\overline{X}_n$.

If $U$ has exactly three neighbors in $P$, by Lemma~\ref{lem2.5}, we
have $|N_{AQ_n} (U)\cap N_{AQ_n}(X,Y,Z)| \leq 8$. Therefore,
$|N_{AQ_n}(U,X,Y,Z)| \geq 6n-17-1+2n-1-3-8=8n-30$, and
$|N_{AQ_n}(U,X,Y,Z)|\geq 6n-15-1+2n-1-3-8=8n-28$ if
$Z=\overline{X}_n$.

The lemma follows.
\end{pf}

\section{Main Results}
\par
In this section, we present our main results, that is, we determine
the $2$-extra connectivity and the $2$-extra edge-connectivity of
the augmented cube $AQ_n$.

\vskip6pt\begin{thm}\
 $\kappa_2(AQ_n)=6n-17$ for $n\geq9$.
\end{thm}

\begin{pf}\
Take a path $P=(\overline{X}_i, X, \overline{X}_{i+2})$ in $AQ_n$,
where $2\leq i\leq {n-3}$. Then, by (\ref{e2.4}) in the proof of
Lemma~\ref{lem2.6}, $|N_{AQ_n}(P)|=6n-17$.

Let $H=AQ_n-(P\cup N_{AQ_n}(P))$. Then, for $n\geq9$,
 $$
 \begin{array}{rl}
 |V(H)|&=|V(AQ_n)|-|V(P)|-|N_{AQ_n}(P)|\\
  &=2^n-3-(6n-17)\\
   &=2^n-6n+14>0,
  \end{array}
 $$
that is, $V(H)\neq \emptyset$. By Lemma~\ref{lem2.5},
$|N_{AQ_n}(Y)\cap N_{AQ_n}(P)|\leq12$ for any $Y\in V(H)$ and
$|N_{AQ_n}(e) \cap N_{AQ_n}(P)|\leq 24$ for any $e\in E(H)$. It
follows that, for $n\geq9$,
 $$
 \begin{array}{rl}
 &|N_{AQ_n}(Y)\cap N_{AQ_n}(P)|\leq12<2n-1=|N_{AQ_n}(Y)|,\ {\rm and}\\
 & |N_{AQ_n}(e) \cap N_{AQ_n}(P)|\leq 24<4n-8\leq|N_{AQ_n}(e)|,
 \end{array}
 $$
which mean that there is neither isolated
vertex nor isolated edge in $AQ_n-N_{AQ_n}(P)$, and so $\kappa_2(AQ_n)\leq6n-17$
for $n\geq9$.

\vskip6pt

Now we only need to prove $\kappa_2(AQ_n)\geq6n-17$ for $n\geq9$.

Suppose that there is a subset $S\subset V(AQ_n) $ with $|S|
\leq6n-18$ such that there is neither isolated vertex nor
isolated edge in $AQ_n-S$. We want to deduce a contradiction by
proving that $AQ_n-S$ is connected.

Let $AQ_n=L\odot R$, where $L=AQ^{0}_{n-1}$ and $R=AQ^{1}_{n-1}$.
For convenience, let $S_L=S\cap L$ and $S_R=S \cap R$. Without loss
of generality we may suppose that $|S_L|\leq|S_R|$. Then
 \begin{equation}\label{e3.1}
 |S_L|\leq(6n-18)/2=3n-9.
 \end{equation}

We prove that $AQ_n-S$ is connected by two steps. In step 1, we
prove that $L-S_L$ is connected in $AQ_n-S$. And in step 2, we prove
that any vertex in $R-S_R$ can be connected to some vertex in
$L-S_L$.

\vskip6pt

{\bf Step 1.}  $L-S_L$ is connected in $AQ_n-S$.

In this case, by (\ref{e3.1}) and Lemma~\ref{lem2.5a}, for $n\geq9$,
 \begin{equation}\label{e3.2}
|S_L| \leq3n-9<4n-12=\kappa_1(L).
 \end{equation}

If there are no isolated vertices in $L-S_L$, then $L-S_L$ is
connected by (\ref{e3.2}).

Suppose now that there exist isolated vertices in $L-S_L$. Note that
$L$ is $(2n-3)$-regular since $L\cong AQ_{n-1}$. By
Lemma~\ref{lem2.5}, any two vertices in $L$ have at most four common
neighbors. To isolate two vertices in $L$, we have to remove at
least $4n-12$ vertices. By (\ref{e3.2}), there is exactly one
isolated vertex, say $X$ in $L-S_L$. Thus
  \begin{equation}\label{e3.3a}
|S_L|\geq|N_L(X)|=2n-3,
 \end{equation}
and so,
 $
 |S_R|=|S|-|S_L| \leq6n-18-(2n-3)=4n-15,
 $
that is,
 \begin{equation}\label{e3.3}
|S_R|\leq 4n-15.
 \end{equation}
Let $S_L' = S_L\cup \{X\}$. Then $L-S'_L$ contains no isolated
vertices. By (\ref{e3.2}), $|S'_L|<\kappa_1(L)$ for $n\geq9$. Thus,
$L-S'_L$ is also connected. We only need to show that $X$ can be
connected to some vertices in $L-S'_L$ via some vertices in $R-S_R$.

By Definition~\ref{Def2.1}, $X$ has two neighbors in $R$, that is,
$X_n$ and $\overline{X}_n$. By our hypothesis there are no isolated
vertices in $AQ_n-S$, then there is at least one in
$\{X_n,\overline{X}_n\}$ is not in $S_R$. Without loss of
generality, assume that $X_n$ is not in $S_R$. Consider two cases
according as the vertex $\overline{X}_n$ is in $S_R$ or not.

\vskip6pt

{\bf Case 1.1.}\ $\overline{X}_n\notin S_R$.

Since $X_n$ and $\overline{X}_n$ are adjacent in $AQ_n$, and
$X\overline{X}_n$ is a complement edge of dimension $n$, by
(\ref{e2.2}), $N_{AQ_n}(X_n)\cap N_{AQ_n}(\overline{X}_n)
=\{X,\overline{X}_{n-1},(X_n)_{n-1},(\overline{X}_{n})_{n-1}\}$, that is, $|N_R(X_n)\cap
N_{R}(\overline{X}_n)|=2$. By (\ref{e3.3}), we have
 $$
 |N_R( X_n, \overline{X}_n)|=2(2n-3)-2-2=4n-10>4n-15\geq |S_R|.
 $$
Thus $|N_R(X_n, \overline{X}_n)\setminus S_R| \geq 5$. Let
 $$
 Y =\{Y^{i}: Y^{i} \in N_R(X_n, \overline{X}_n)\setminus S_R \}.
 $$
Then $|Y|\geq 5$. If there is some vertex $Y^{i}\in Y$ such that at
least one of $Y^{i}_{n}$ and $\overline{Y^{i}}_{n}$ is not $S_L$,
then we are done. So assume that both $Y^{i}_{n}$ and
$\overline{Y^{i}}_{n}$ are in $S_L$ for any $Y^{i}\in Y$.
Take $Y^1, Y^2\in Y$ such that $(X_n, \overline{X}_n, Y^{1})$ (or
$(\overline{X}_n, X_n, Y^{1})$) is a path and $Y^2$ is adjacent to
one in $\{X_n, \overline{X}_n, Y^{1}\}$. By Lemma~\ref{lem2.7}, we
have that
 $$
 |N_R(X_n, \overline{X}_n, Y^{1}, Y^{2})| \geq 8(n-1)-29=8n-37.
 $$
Let
 $$
 C = N_R(X_n, \overline{X}_n, Y^{1}, Y^{2})\cup \{\overline{X}_n, Y^{1}, Y^{2} \}.
 $$
Then $|C|\geq 8n-34$. Let $E^h_{n}=\{ UU_n: U\in C \}$. Noting that
all edges in $E^h_n$ are hypercube edges of dimension $n$, we have
that, for $n\geq 9$,
 $$
 |E^h_{n}| = |C| \geq8n-34>6n-18\geq |S|,
 $$
which means that there exists an edge, say $UU_n$, in $E^h_{n}$ such
that neither of its two end-vertices is in $S$. Since $X_n$,
$\overline{X}_n$, $Y^{1}$ and $Y^{2}$ are all not in S, $X$ can be
connected to $L-S_L^{'}$ via vertices in $R-S_R$ and the edge
$UU_n$.

\vskip6pt

{\bf Case 1.2.}\  $\overline{X}_n\in S_R$.

Since there are no isolated edges in $AQ_n-S$, we have
$N_R(X_n)\setminus S_R\neq \emptyset$. If there is some $U\in
N_R(X_n)\setminus S_R$ such that at least one in $\{U_n,
\overline{U}_n\}$ is not in $S_L$, then we are done. So assume that
$U_n$ and $\overline{U}_n$ are both in $S_L$ for any $U\in
N_R(X_n)\setminus S_R$. Noting that $\kappa_1(R)=4(n-1)-8$ and
(\ref{e3.3}), we have that
 $$
 |N_R(X_n, U)|\geq4n-12>4n-15\geq |S_R|,
 $$
which implies $N_R(X_n, U)\setminus S_R \neq \emptyset$. Let
 $$
 Z=\{Z^{i}: Z^{i}\in N_R(X_n, U)\setminus S_R \}.
 $$
Then $|Z|\geq 2$.
Take $Z^1, Z^2\in Z$ such that $(X_n, U, Z^{1})$ (or $(U, X_n,
Z^{1})$) is a path and $Z^2$ is adjacent to one in $\{X_n, U,
Y^{1}\}$. By Lemma~\ref{lem2.7}, we have that
 $$
 |N_R(X_n, U, Z^{1},Z^{2})|\geq 8(n-1)-31=8n-39.
 $$
Let
 $$
 D=N_R(X_n, U, Z^{1}, Z^{2})\cup \{ X_n, U, Z^{1}, Z^{2}\}.
 $$
Then $|D| \geq8n-35$. Let $E^{c}_{n}=\{A\bar{A}_n: A\in D\}$. Noting
that all edges in $E^c_n$ are complement edges of dimension $n$, we
have that, for $n\geq9$,
 $$
 |E^{c}_{n}|=|D|\geq8n-35>6n-18.
 $$
Hence, there exists an edge, say $A\bar{A}_n$, in $E^{c}_{n}$ such
that neither of its two end-vertices is $S$. Since $X_n$, $U$,
$Z^{1}$ and $Z^{2}$ are all not in $S$, $X$ can be connected to
$L-S_L^{'}$ via vertices in $R-S_R$ and the edge $A\bar{A}_n$.

\vskip6pt

{\bf Step 2.}\ Any vertex in $R-S_R$ can be connected to some vertex
in $L-S_L$.

Let $U$ be any vertex in $R-S_R$. Consider $U_n$ and
$\overline{U}_n$, which are neighbors of $U$ in $L$. If at least one
of $U_n$ and $\overline{U}_n$ is not in $S_L$, we are done. So
suppose that both $U_n$ and $\overline{U}_n$ are in $S_L$. Consider
the neighbor $\overline{U}_{n-1}$ of $U$ in $R$. There are two cases
according as $\overline{U}_{n-1}$ is in $S_R$ or not.

\vskip6pt

{\bf Case 2.1.} $\overline{U}_{n-1}\notin S_R$.

Note that $U\overline{U}_{n-1}$ is a complement edge of dimension
$(n-1)$ in $AQ_n$ and
 $$
 (N_{AQ_n}(U)\cap N_{AQ_n}(\overline{U}_{n-1}))\cap L=\{U_n, \overline{U}_n\}\subseteq S_L.
 $$
Since $U\overline{U}_{n-1}$ is not an isolated edge in $AQ_n-S$,
there exists a vertex $V\in N_R(U, \overline{U}_{n-1})\setminus
S_R$. Then $V_n$ and $\overline{V}_{n}$ are neighbors of $V$ in $L$.
If at least one of $V_n$ and $\overline{V}_{n}$ is not in $S_L$, we
are done. So assume that both $V_n$ and $\overline{V}_{n}$ are in
$S_L$. By (\ref{e2.4}) and (\ref{e2.5}) in the proof of Lemma~\ref{lem2.6}, we have
that
 \begin{equation}\label{e3.4}
 |N_R(U, \overline{U}_{n-1}, V)|\geq6(n-1)-15 =6n-21.
 \end{equation}
Since
 $$
 |S_L|\geq |\{U_n, \overline{U}_{n}, V_n,\overline{V}_{n}\}|=4,
 $$
we have that
  \begin{equation}\label{e3.5}
  |S_R|=|S|-|S_L|\leq6n-18-4=6n-22,
  \end{equation}
By (\ref{e3.4}) and (\ref{e3.5}), we have
 $
 |N_R(U,\overline{U}_{n-1}, V) |> |S_R|,
 $
that is,
 $$
 N_R (U, \overline{U}_{n-1}, V)\setminus S_R \neq \emptyset.
 $$
If there is some $W\in N_R(U, \overline{U}_{n-1}, V)\setminus S_R$
such that at least one of $W_n$ and $\overline{W}_{n}$, which are
two neighbors of $W$ in $L$, is not in $S_L$, we are done. So assume
that both $W_n$ and $\overline{W}_{n}$ are in $S_L$ for any $W \in
N_R(U, \overline{U}_{n-1}, V)\setminus S_R$. By Lemma~\ref{lem2.7},
we have that
 $$
 |N_R(U, \overline{U}_{n-1}, V, W)| \geq 8n-37.
 $$
Let
 $$
 C'=N_R(U, \overline{U}_{n-1}, V, W) \cup \{U, \overline{U}_{n-1}, V, W\}.
 $$
Then $|C'| \geq8n-33$. Let $E'_{n}=\{AA_n: A\in C'\}$. Noting that
all edges in $E'_n$ are hypercube edges of dimension $n$, we have
that, for $n\geq 9$,
 $$
 |E'_{n}|=|C'| \geq 8n-33>6n-18.
 $$
There exists at least one edge, say $AA_n$, of $E'_{n}$ whose two
end-vertices both are not in $S$. Since $U$, $\overline{U}_{n-1}$,
$V$ and $W$ are all not in $S$, this implies that $U$ can be
connected to $L-S_L$.

\vskip6pt

{\bf Case 2.2.}  $\overline{U}_{n-1}\in S_R$.

Since $U$ is not an isolated vertex in $AQ_n-S$ and two neighbors
$U_n$ and $\overline{U}_n$ of $U$ in $L$ are both in $S$, there
exists a vertex $B \in N_R(U)\setminus S_R$. Then $B_n$ and
$\overline{B}_{n}$ are neighbors of $B$ in $L$. If at least one of
$B_n$ and $\overline{B}_{n}$ is not in $S_L$, we are done. So assume
that both $B_n$ and $\overline{B}_{n}$ are in $S_L$. If
$\overline{B}_{n-1}\notin S_R$, we can obtain a path joining $B$ to
some vertex in $L-S_L$ by Case 2.1 by replacing $U$ by $B$.
Therefore assume $\overline{B}_{n-1}\in S_R$ below.

Since $UB$ is not an isolated edge in $AQ_n-S$, there exists a
vertex $F \in N_R(U, B)\setminus S_R$. Then $F_n$ and
$\overline{F}_{n}$ are two neighbors of $F$ in $L$. If at least one
of $F_n$ and $\overline{F}_{n}$ is not in $S_L$, we are done. So
suppose that both $F_n$ and $\overline{F}_{n}$ are in $S_L$. By
Lemma~\ref{lem2.6}, we have that
 \begin{equation}\label{e3.6}
 |N_R(U, B, F)| \geq6(n-1)-17=6n-23.
 \end{equation}
Since
 $$
 |S_L| \geq |\{U_n, \overline{U}_n, B_n, \overline{B}_{n},
 F_n, \overline{F}_{n}\}|=6,
 $$
we have that
 \begin{equation}\label{e3.7}
 |S_R|=|S|-|S_L|\leq6n-18-6=6n-24.
 \end{equation}
Comparing (\ref{e3.6}) with (\ref{e3.7}), we have that $N_R(U, B,
F)\setminus S_R \neq \emptyset$. Let $Q\in N_R(U, B,
F)\setminus S_R$. By Lemma~\ref{lem2.7}, we have that
 $$
 |N_R(U, B, F, Q)| \geq8n-39.
 $$
Let
 $$
 C''=N_R(U, B, F, Q) \cup\{U, B, F, Q\}.
 $$
Then $|C''| \geq8n-35$. Let $E''_{n}=\{AA_n: A\in C''\}$. Noting
that all edges in $E''_n$ are hypercube edges of dimension $n$, we
have, for $n\geq9$,
 $$
 |E''_{n}|= |C''| \geq8n-35>6n-18.
 $$
There exists an edge, say $AA_n$, of $E''_{n}$ whose two
end-vertices both are not in $S$. Since $U$, $B$, $F$ and $Q$ are
all not in $S$, thus $U$ can be connected to $L-S_L$.

The proof of the theorem is complete.
\end{pf}

\vskip6pt\begin{thm}\
 $\lambda_2(AQ_n)=6n-9$ for $n\geq4$.
\end{thm}

\begin{pf}\
Let $C_3$ be a cycle of length three in $AQ_n$, $U$ be any vertex
not in $C_3$, and let $e$ be any edge $e$ not incident with any
vertex in $C_3$. Obviously, any vertex not in $C_3$ can have at most
3 neighbors in $C_3$. Thus, for $n\geq4$,
 $$
 |E_{AQ_n}(U) \cap E_{AQ_n}(C_3)|\leq 3<2n-1=|E_{AQ_n}(U)|,
 $$
and
 $$
 |E_{AQ_n}(e) \cap E_{AQ_n}(C_3)|\leq 6<4n-4=|E_{AQ_n}(e)|.
 $$
So, there are no isolated vertices or isolated edges in
$AQ_n-E_{AQ_n}(C_3)$. That is, $E_{AQ_n}(C_3)$ is a $2$-edge-cut of
$G$. It follows that, for $n\geq4$,
 $$
 \lambda_2(AQ_n)\leq E_{AQ_n}(C_3)=6n-9.
 $$

In the following, we only need to prove that $\lambda_2(AQ_n)\geq
6n-9$ for $n\geq4$.

Let $F$ be an arbitrary $2$-edge-cut in $AQ_n$ with $|F|\leq 6n-10$
such that there are neither isolated vertices nor isolated edges in
$AQ_n-F$. Let $F_L=F\cap L$ and $F_R=F\cap R$. Without loss of
generality we may suppose that $|F_L|\leq |F_R|$. Then
 $$
 |F_L|\leq\frac 12\,(6n-10)=3n-5.
 $$
We will deduce a contradiction by proving that $AQ_n-F$ is connected
by two steps. In step 1, we show that $L-F_L$ is connected in
$AQ_n-F$. In step 2, we show that any vertex of $R$ can be connected
to $L$ in $AQ_n-F$.

\vskip6pt

{\bf Step 1.}  $L-F_L$ is connected in $AQ_n-F$.

By our hypothesis and Lemma~\ref{lem2.5a}, for $n\geq4$, we have that
 \begin{equation}\label{e3.8}
 |F_L|\leq 3n-5< 4(n-1)-4=\lambda_1(L).
 \end{equation}
Thus, if there are no isolated vertices in $L-F_L$, then $L-F_L$ is
connected, and so we are done. In the following discussion, we
assume that there exists an isolated vertex $X$ in $L-F_L$.

Since $L$ is $(2n-3)$-regular and any two vertices are incident with
at most one edge, to get two isolated vertices in $L$, we have to
remove at least $4n-7$ edges from $L$. However, by (\ref{e3.8}),
$|F_L|\leq3n-5<4n-7$ for $n\geq 4$. This shows that there is just
one isolated vertex $X$ in $L-F_L$. Then by Lemma~\ref{lem2.3}, we have
  $$
  \lambda(L-X)\geq \kappa(L-X) \geq \kappa(L)-1=\left\{ \begin{array}{lll}
     4-1=3 &\ {\rm if }\   n=4\\
     2n-3-1=2n-4  &\ {\rm if }\  n>4\\
     \end{array} \right.
 $$
and
 $$
 |F_L|-|E_L(X)|\leq3n-5-(2n-3)=n-2,
 $$
which implies that $|F_L|-|E_L(X)|< \lambda(L-X)$ for $n\geq 4$. In other
words, the subgraph $H=(L-X)-F_L=(L-F_L)-X$ is connected. In the
following we only need to prove that $X$ can be connected to $H$ in
$AQ_n-F$. Since $AQ_n-F$ contains no isolated vertices, at least one
of two edges $XX_n$ and $X\overline{X}_n$ is not in $F$. Without
loss of generality, we may assume that $XX_n$ is not $F$. Consider
two cases according as the edge $X\overline{X}_n$ is in $F$ or not.

\vskip6pt

{\bf Case 1.1.}\ $X\overline{X}_n \notin F$.

Note that $X_n\overline{X}_{n-1}$ and
$\overline{X}_{n}\overline{X}_{n-1}$ are edges in $AQ_n$, where
$\overline{X}_{n-1}$ is in $L$. If at least one of the two edges is
not in $F$, we are done. So, we can assume that both the two edges
are in $F$. We will construct $4n-8$ edge disjoint paths joining $X$
to some vertex in $L-X$ in the following.


 Let
 $$
 \begin{array}{rl}
 &E_1=\{X_nX^i:\ X_nX^i \in
 E_R(X_n)\backslash\{X_n\overline{X}_{n}\}\},\\
 & E_2=\{\overline{X}_{n}Y^j:\ \overline{X}_{n}Y^j \in
 E_R(\overline{X}_{n})\backslash \{X_n\overline{X}_{n}\}\},\\
 & F'=F\setminus (E_L(X) \cup\{X_n\overline{X}_{n-1},
 \overline{X}_{n}\overline{X}_{n-1}\}).
 \end{array}
 $$
Then $|E_1|=|E_2|=2n-4$, $E_1 \cap E_2=\emptyset$, and
 \begin{equation}\label{e3.9}
 |F'|\leq 6n-10-(2n-3)-2=4n-9.
 \end{equation}
Let
 \begin{equation}\label{e3.10}
 P_i=(X, X_n, X^i, X^i_n)\ \ {\rm and}\ Q_j=(X,
 \overline{X}_{n}, Y^j, \overline{Y^j}_n)
 \end{equation}
be a path joining $X$ to some vertex in $L-X$, and let
 $$
 \mathscr P=\{P_i: 1\leq i\leq 2n-4\}\cup\{Q_j: 1\leq j\leq 2n-4\}.
 $$
Then
 \begin{equation}\label{e3.11}
 |\mathscr P|=4n-8.
 \end{equation}
Since these paths defined in (\ref{e3.10}) are edges disjoint,
comparing (\ref{e3.9}) and (\ref{e3.11}), we have that there exists
a path $P\in \mathscr P$ such that $E(P) \cap F'=\emptyset$. Then
$X$ can be connected to a vertex in $L$.

\vskip6pt

{\bf Case 1.2.}\ $X\overline{X}_n \in F$.

If $X_n\overline{X}_{n-1} \notin F$, we are done. So assume that
$X_n\overline{X}_{n-1} \in F$. Since there are no isolated edges in
$AQ_n-F$, we have $E_R(X_n)-F_R \neq \emptyset$.

If $X_n\overline{X}_{n} \notin F_R$, we can obtain a path joining
$X$ via a path to a vertex in $L-X$ by the $4n-8$ paths $\mathscr P$
constructed in Case 1.1. Hence, we can assume $X_n\overline{X}_{n}
\in F_R$ below. We will construct $4n-9$ edge disjoint paths joining
$X$ to some vertex in $L-X$ in the following.

Let $X_nW\in E_R(X_n)-F_R$. Then $W\neq \overline{X}_{n}$.

Let
$$
 \begin{array}{rl}
 &E'_1=\{WX^i:\ WX^i \in E_R(W)\backslash\{X_nW\}\},\\
 & E'_2=\{X_nY^j:\ X_nY^j \in E_R(X_n)\backslash\{X_nW, X_n\overline{X}_{n}\}\},\\
 & F^{*}=F\setminus (E_L(X) \cup\{ X\overline{X}_{n}, X_n\overline{X}_{n-1},
 X_n\overline{X}_{n} \}).
 \end{array}
 $$
Then $|E'_1|=2n-4$, $|E'_2|=2n-5$ and
\begin{equation}\label{e3.12}
 |F^{*}| \leq 6n-10-(2n-3)-3=4n-10.
 \end{equation}
Let
 \begin{equation}\label{e3.13}
 P^{*}_i=(X, X_n, W, X^i, X^i_n)\ \ {\rm and}\ Q^{*}_j=(X, X_n,
Y^{i}, \overline{Y^j}_n))
 \end{equation}
be a path joining $X$ to some vertex in $L-X$, and let
$$
 \mathscr {P^{*}}=\{P^{*}_i: 1\leq i\leq2n-4\}\cup\{Q^{*}_j: 1\leq j\leq2n-5\}.
 $$
Then
 \begin{equation}\label{e3.14}
 |\mathscr {P^{*}}|=4n-9.
 \end{equation}
Since these paths defined in (\ref{e3.13}) are edges disjoint,
comparing (\ref{e3.12}) and (\ref{e3.14}), there exists a path
$P^{*}\in \mathscr {P^{*}}$ such that $E(P^{*}) \cap
F^{*}=\emptyset$. This implies that vertex $X$ can be connected to a
vertex in $L-X$.

\vskip6pt

{\bf Step 2.} Any vertex $X$ of $R$ can be connected to $L$ in
$AQ_n-F$.

Suppose that $X$ is an arbitrary vertex in $R$, if $\{XX_n,
X\overline{X}_{n}\}\nsubseteq F$, where $X_n$ and $\overline{X}_{n}$
are both in $L$, we are done. Thus, assume that $\{XX_n,
X\overline{X}_{n}\}\subseteq F$. Since there is neither isolated
vertex nor isolated edge in $AQ_n-F$, the vertex $X$ lies on a path
$T$ of length $2$ in $R-F_R$. Assume $V(T)=\{X,Y,Z\}$. If $\{YY_n, Y\overline{Y}_{n}, ZZ_n,
Z\overline{Z}_{n}\}\nsubseteq F$, we are done. Hence, assume that
$\{YY_n, Y\overline{Y}_{n}, ZZ_n, Z\overline{Z}_{n}\}\subseteq F$ in
the following.

Let $F^{*}=F\setminus \{XX_n, X\overline{X}_{n}, YY_n,
Y\overline{Y}_{n}, ZZ_n, Z\overline{Z}_{n}\}$.

Then
 \begin{equation}\label{e3.16}
 |F^*|\leq 6n-10-6=6n-16.
 \end{equation}
Note that
\begin{equation}\label{e3.17}
|E_R(T)|\geq 3*(2n-5)=6n-15.
 \end{equation}

By Lemma~\ref{lem2.6}, we have $|N_R (X) \cap N_R (Y) \cap N_R
(Z)|\leq1$. We will construct edge disjoint paths joining $V(T)$ to
$V(L)$ according to the following two cases.

\vskip6pt

{\bf Case 2.1} $|N_R (X) \cap N_R (Y) \cap N_R (Z)|= 0$.

In this case, each vertex in $N_R(T)$ is incident to at most two
edges in $E_R(T)$. Since every vertex is incident to exactly two
crossed edges in $AQ_n$, we can construct $|E_R(T)|$ edge disjoint
paths $\mathscr P$ joining $V(T)$ to $V(L)$ as follows.

For any vertex $W^i$ in $N_R(T)$, if $W^i$ is incident to
exactly one edge $AW^i$ in $E_R(T)$ where $A\in V(T)$, let
$P_i=(A, W^i, W^i_n)$; if $W^i$ is incident to exactly two edges $AW^i$
and $BW^i$ in $E_R(T)$ where $A,B\in V(T)$, let $P_i=(A, W^i, W^i_n)$ and
$P_i'=(B, W^i, \overline {W^i}_n)$.

Comparing (\ref{e3.16}) and (\ref{e3.17}), there exists a path $P\in
\mathscr {P}$ such that $E(P) \cap F^{*}=\emptyset$. This implies
that vertex $X$ can be connected to a vertex in $L$.

\vskip6pt

{\bf Case 2.2} $|N_R (X) \cap N_R (Y) \cap N_R (Z)|= 1$.

In this case, there is a vertex $U$ in $N_R(T)$ incident to exactly
three edges in $E_R(T)$ and other vertices in $N_R(T)$ is incident
to at most two edges in $E_R(T)$.

Since $|F_R|\leq |F^*|$, comparing (\ref{e3.16}) and (\ref{e3.17}),
there exists an edge $e\in E_R(T)\setminus F_R$. Without loss of
generality, we may assume $e=XW$. If $\{WW_n,
W\overline{W}_{n}\}\nsubseteq F$, where $W_n$ and $\overline{W}_{n}$
are both in $L$, we are done. Thus, assume that $\{WW_n,
W\overline{W}_{n}\}\subseteq F$.

Let $F'=F\setminus \{XX_n, X\overline{X}_{n}, YY_n,
Y\overline{Y}_{n}, ZZ_n, Z\overline{Z}_{n}, WW_n,
W\overline{W}_{n}\}$.

Then
 \begin{equation}\label{e3.18}
 |F'|\leq 6n-10-8=6n-18.
 \end{equation}

Since every vertex is incident to exactly two crossed edges in
$AQ_n$, we can construct edge disjoint paths
$\mathscr P$ joining $V(T)$ to $V(L)$ as follows.

For any vertex $W^i$ in $N_R(T)\setminus \{W\}$, if $W^i$ is
incident to exactly one edge $AW^i$ in $E_R(T)$ where $A\in V(T)$,
let $P_i=(A, W^i, W^i_n)$; if $W^i$ is incident to exactly two edges
$AW^i$ and $BW^i$ in $E_R(T)$ where $A,B\in V(T)$, let
$P_i=(A, W^i, W^i_n)$ and $P_i'=(B, W^i, \overline {W^i}_n)$; if $W^i=U$, let
$P_i=(Y, U, U_n)$ and $P_i'=(Z, U, \overline {U}_n)$.

\vskip6pt

{\bf Case 2.2.1} The vertex $W$ is incident to one edge in $E_R(T)$. Hence, we
construct at least $|E_R(T)|-2\geq6n-17$ edge disjoint paths $\mathscr P$ jointing
$V(T)$ to $V(L)$. By (\ref{e3.18}), there exists a path $P\in
\mathscr {P}$ such that $E(P) \cap F'=\emptyset$. This implies
that vertex $X$ can be connected to a vertex in $L$.

\vskip6pt

{\bf Case 2.2.2} The vertex $W$ is incident to at least two edges
in $E_R(T)$. Hence, we construct at least $|E_R(T)|-3\geq6n-18$
edge disjoint paths $\mathscr P$ jointing $V(T)$ to $V(L)$.
If there exists a path  $P\in
\mathscr {P}$ such that $E(P) \cap F'=\emptyset$, we are done.
Assume $E(P) \cap F'\neq\emptyset$ for every path in $\mathscr P$.
By (\ref{e3.18}), the faulty edges $F'$ are all in $\mathscr P$.
Let $P=(X, W, W_j, (W_j)_n)$, where $W_j$ is not in $T$ and $WW_j\nsubseteq F$.
Then, $P$ is fault-free. This implies that vertex $X$ can be connected to
a vertex in $L$.

\par
We proved that $AQ_n-F$ is connected, which means
$\lambda_2(AQ_n)\geq6n-9$ for $n\geq4$.

The theorem follows.
\end{pf}

\section{Conclusions}

\par
In this paper, we explore two stronger measurement parameters for
the reliability and the tolerance of networks called the 2-extra
connectivity $\kappa_2(G)$ and the 2-extra edge-connectivity
$\lambda_2(G)$ of a connected graph $G$, which not only compensate
for some shortcomings but also generalize the classical connectivity
$\kappa(G)$ and the classical edge-connectivity $\lambda(G)$, and so
can provide more accurate measures for the reliability and the
tolerance of a large-scale parallel processing system. The augmented
cube $AQ_n$, as an important variant of the hypercube $Q_n$, has
many desirable properties (for more results, see, for example,
\cite{xml07, q16, q17, cs00, cs01, cs02, cs08, f10, hc10, kth10,
mlx08, mtxl09, xwm07}). Here, we have showed that
$\kappa_2(AQ_n)=6n-17$ for $n\geq9$; and $\lambda_2(AQ_n)=6n-9$ for
$n\geq4$. In other words, for $n\geq9$ (respectively, $n\geq 4$), at
least $6n-17$ vertices (respectively, $6n-9$ edges) of $AQ_n$ have
to be removed to get a disconnected graph that contains no isolated
vertices and isolated edges. Compared with previous results, our
results enhance the fault tolerant ability of this kind of network
theoretically. 　　


\begin{thebibliography}{10}


\bibitem{bbst81}
D. Bauer, F. Boesch, C. Suffel and R. Tindell, Connectivity extremal
problems and the design of reliable probabilistic networks, The
theory and Application of Graphs, New York: Wiley (1981) 89-98.

\bibitem{q3}
F. T. Boesch, Synthesis of reliable networks-a survey, IEEE
Transactions on Reliability. 35 (1986) 240-246.

\bibitem{q16}
M. Chan, The distinguishing number of the augmented cube and
hypercube powers, Discrete Mathematics. 308 (2008) 2330-2336.


\bibitem{q17}
Y. C. Chen, M. H. Chen, J. J. M. Tan, Maximally local connectivity
on augmented cubes, Lecture Notes in Computer Science. 5574 (2009)
121-128.

\bibitem{cls10}
E. Cheng, L. Lipt\'ak, F. Sala, Linearly many faults in
2-tree-generated networks, Networks. 55 (2) (2010) 90-98.


\bibitem{cs00}
S. A. Choudum, V. Sunitha, Wide-diameter of augmented cubes,
Technical report, Department of Mathematics, Indian Institute of
Technology Madras, Chennai, November 2000.

\bibitem{cs01}
S. A. Choudum, V. Sunitha, Automorphisms of augmented cubes,
Technical report, Department of Mathematics, Indian Institute of
Technology Madras, Chennai, March 2001.


\bibitem{cs02}
S. A. Choudum, V. Sunitha, Augmented cubes, Networks. 40 (2) (2002)
71-84.

\bibitem{cs08}
S. A. Choudum and V. Sunitha, Automorphisms of augmented cubes. Int.
J. Comput. Math., 85 (11) (2008), 1621-1627.

\bibitem{es89}
A. H. Esfahanian, Generalized measures of fault tolerance with
application to $n$-cube networks, IEEE Transactions on computers. 38
(11) (1989) 1586-1591.

\bibitem{q5}
A. H. Esfahanian, S. L. Hakimi, On computing a conditional
edge-connectivity of a graph, Information Processing Letters. 27
(1988) 195-199.

\bibitem{q7}
J. F\`{a}brega, M. A. Fiol, Extraconnectivity of graphs with large
girth, Discrete Mathematics. 127 (1994) 163-170.

\bibitem{f10}
J.-S. Fu, Edge-fault-tolerant vertex-pancyclicity of augmented
cubes. Information Processing Letters, 110 (11) (2010), 439-443.


\bibitem{q12}
A. Hellwig, D. Rautenbach, L. Volkmann, Note on the connectivity of
line graphs, Information Processing Letters. 91 (2004) 7-10.


\bibitem{q11}
A. Hellwig, L. Volkmann, Sufficient conditions for
$\lambda'$-optimality in graphs of diameters 2, Discrete
Mathematics. 283 (2004) 113-120.


\bibitem{hc10}
S.-Y., Hsieh, and Y.-R. Cian, Conditional edge-fault Hamiltonicity
of augmented cubes. Information Sciences, 180 (13) (2010),
2596-2617,



\bibitem{hy97}
S.-C. Hu, C.-B. Yang, Fault tolerance on star graphs, International
Journal of Foundations of Computer Science. 8(2) (1997) 127-142.

\bibitem{kth10}
T.-L. Kung, Y.-H. Teng, and L.-H. Hsu, The panpositionable
panconnectedness of augmented cubes. Information Sciences, 180 (19)
(2010), 3781-3793



\bibitem{mlx08}
M. Ma, G. Liu, J.-M. Xu, The super connectivity of augmented cubes,
Information Processing Letters. 106(2) (2008) 59-63.

\bibitem{mtxl09}
M. Ma, X. Tan, J.-M. Xu, G. Liu, A note on ``The super connectivity
of augmented cubes", Information Processing Letters. 109(12) (2009)
592-593.

\bibitem{wz09}M. Wan, Z. Zhang, A kind of conditional vertex connectivity of
star graphs, Applied Mathematics  Letters. 22 (2009) 264-267.

\bibitem{q2}
J.-M. Xu, Theory and Application of Graphs, Kluwer Academic
Publishers, Dordrecht Boston London. 2003.


\bibitem{xlmh05}
J.-M. Xu, M. L\"u, M. Ma, A. Hellwig, Super connectivity of line
graphs, Information Processing Letters. 94 (4) (2005) 191-195.

\bibitem{xml07}
J.-M. Xu, M. Ma, G. Liu, Panconnectiviy and edge-fault-tolerant
pancyclicity of augmented cubes, Parallel Computing. 33 (2007)
36-42.

\bibitem{xwm07}
J.-M. Xu, W.-W.Wang, M.-J. Ma, Fault-tolerant pancyclicity of
augmented cubes, Information processing Letters. 103 (2007) 52-56.

\bibitem{xww10}
J.-M. Xu, J.-W. Wang, W.-W. Wang, Super and restricted connectivity
of some interconnection networks, Ars Combinatoria. 94 (2010) 25-32.

\bibitem{q13}
J.-M. Xu, K.-L. Xu, On restricted edge-connectivity of graphs,
Discrete Mathematics. 243 (1-3) (2002) 291-298.


\bibitem{xzhz05}
J.-M. Xu,  Q. Zhu, X. Hou and T. Zhou, On restricted connectivity
and extra connectivity of hypercubes and folded hypercubes, Journal
of Shanghai Jiaotong University. E-10 (2) (2005) 203-207.

\bibitem{ylg10}
W.-H. Yang, H.-Z. Li, X.-F. Guo, A kind of conditional fault
tolerance of $(n, k)$-star graph, Information Processing Letters.
110 (22) (2010) 1007-1011.


\bibitem{zxy10}
Z. Zhang, W. Xiong, W.-H. Yang, A kind of conditional fault
tolerance alternating group graphs, Information Processing Letters.
110 (2010) 998-1002.

\bibitem{zx06}
Q. Zhu, J.-M. Xu, On restricted edge-connectivity and extra
edge-connectivity of hypercubes and folded hypercubes, Journal of
University of Science and Technology of China. 36 (3) (2006)
249-253.

\bibitem{zxhx07}
Q. Zhu, J.-M. Xu, X. Hou, M. Xu, On reliability of the folded
hypercubes, Information Sciences. 177 (8) (2007) 1782-1788.



\end{thebibliography}
\end{document}